\documentclass[11pt]{amsart}
\usepackage{amssymb,latexsym}
\usepackage{graphics}

\newcommand{\s}{\mathfrak{S}}
\newcommand{\I}{\mathcal{I}}
\newcommand{\hh}{\mathsf{h}}
\newcommand{\rr}{\mathsf{r}}
\newcommand{\dd}{\mathsf{d}}
\newcommand{\uu}{\mathsf{u}}
\newcommand{\paths}{\mathsf{Sh}}

\numberwithin{equation}{section}
\newtheorem{theorem}{Theorem}[section]

\newtheorem{corollary}[theorem]{Corollary}

\newtheorem{question}[theorem]{Open Problem}

\newtheorem{lemma}[theorem]{Lemma}

\def\ls{\leq}

\def\SS{\s}

\title[Symmetric Schr\"oder paths and restricted involutions]
{
Symmetric Schr\"oder paths and restricted involutions
}
\author{Eva Y. P. Deng, Mark Dukes, Toufik Mansour and Susan Y. J. Wu }
\address{Department of Applied Mathematics, Dalian University of Technology, Dalian 116024, P.R. China.}
\email{ypdeng@dlut.edu.cn}
\address{Science Institute, University of Iceland, Reykjav\'{i}k, Iceland.}
\email{dukes@hi.is}
\address{Department of Mathematics, University of Haifa, Haifa 31905, Israel.}
\email{toufik@math.haifa.ac.il}
\address{Department of Mathematics, Tianjin Normal University, Tianjin 300387, P.R. China.}
\email{wuyijun@gmail.com}


\begin{document}
\begin{abstract}
Let $\mathcal{A}_k$ be the set of permutations in the symmetric group $\SS_k$ with prefix 12.
This paper concerns the enumeration
of involutions which avoid the set of patterns $\mathcal{A}_k$.
We present a bijection between symmetric Schr\"{o}der paths of length $2n$ and involutions of length $n+1$ avoiding $\mathcal{A}_4$.
Statistics such as the number of right-to-left maxima and fixed points of the involution correspond to the number of steps in the symmetric Schr\"oder path of a particular type.
For each $k\geq 3$ we determine the generating function for the number of involutions avoiding the
subsequences in $\mathcal{A}_k$, according to length, first entry and number of fixed
points.
\end{abstract}


\maketitle
\section{Introduction}
Let $\SS_n$ denote the set of permutations of $[n]=\{1,\ldots,n\}$, written in one-line
notation. For two permutations $\pi\in\SS_n$ and $\tau\in\SS_k$, an {\it occurrence} of
$\tau$ in $\pi$ is a subsequence $1\ls i_1<i_2<\dots<i_k\ls n$ such that $(\pi_{i_1},
\dots,\pi_{i_k})$ is order-isomorphic to $\tau$; in such a context $\tau$ is usually
called a {\it pattern}. We say that $\pi$ {\it avoids} $\tau$, or is $\tau$-{\it
avoiding}, if there is no occurrence of $\tau$ in $\pi$. A natural generalization of
single pattern avoidance is {\it subset avoidance}, that is, we say that $\pi\in\SS_n$
avoids a subset $T\subseteq\SS_k$ if $\pi$ avoids all $\tau\in T$. The set of all
$\tau$-avoiding (resp. $T$-avoiding) permutations of length $n$ is denoted $\SS_n(\tau)$
(resp. $\SS_n(T)$).

Several authors have considered the case of general $k$ in which $T$ enjoys various
algebraic properties. Barcucci et al.~\cite{BDPP} treat the case of permutations avoiding
the collection of permutations in $\SS_k$ that have suffix $(k-1)k$.
Adin and Roichman~\cite{AR} look at the case where $T$ is a Kazhdan--Lusztig cell of
$\SS_k$, or, equivalently, a Knuth equivalence class (see \cite[Vol.~2, Ch.~A1]{St}).
Mansour and Vainshtein~\cite{MV} consider the situation where $T$ is a maximal parabolic
subgroup of $\SS_k$.

In the current paper an analogous result is established for pattern-avoiding involutions.
We say $\pi$ is an {\em involution} whenever $\pi_{\pi_i} = i$ for all $i\in[n]$. Let
$\mathcal{I}_n$ denote the set of involutions of $[n]$. The set of all $\tau$-avoiding
(resp. $T$-avoiding) involutions of length $n$ is denoted $\mathcal{I}_n(\tau)$ (resp.
$\mathcal{I}_n(T)$).

Simion and Schmidt~\cite{SS} considered the first cases of pattern-avoiding involutions,
which was continued in Gouyou-Beauchamps~\cite{G-B} and Gessel~\cite{Ge} for increasing
patterns, and subsequently in Guibert's Ph.D. thesis~\cite{Gu}. This paper concerns 
the enumeration of involutions which avoid the class of
permutations in $\SS_k$ with prefix 12, that is,
$$\mathcal{A}_k=\{\pi_1\pi_2\ldots\pi_k\in \SS_k\mid\pi_1=1,\pi_2=2\}.$$

We denote by $I_n$ the cardinality of the set $\mathcal{I}_n$. We say that $i$ is a fixed
point of a permutation $\pi$ if $\pi_i=i$. Define $J_n(p)$ to be the polynomial
$\sum_{j=0}^n I_{n;j}p^j$, where $I_{n;j}$ is the number of involutions in
$\mathcal{I}_n$ with $j$ fixed points. For example $J_3(p)=3p+p^3$. 
It is not hard to see that the polynomials $J_n(p)$
satisfy the recurrence relation $J_n(p)=pJ_{n-1}(p)+(n-1)J_{n-2}(p)$, $n\geq2$, with
the initial conditions $J_0(p)=1$ and $J_1(p)=p$. The exponential generating function for the sequence
$\{J_n(p)\}_{n\geq 0}$ is given by $e^{px+x^2/2}$.

The main result of this paper can be formulated as follows.

\begin{theorem}\label{thmain}
Let $k\geq2$. The generating function for the number of $\mathcal{A}_k$-avoiding
involutions of length $n$ 
is given by\\
\begin{eqnarray*}
\lefteqn{\displaystyle\sum\limits_{n\geq0}\sum\limits_{\pi\in\I_n (\mathcal{A}_k)}x^np^{\#\mbox{\scriptsize{fixed points}}(\pi)} \;=\;}\\[8pt]
&&\displaystyle\sum_{j=0}^{k-3}J_j(p)x^j-\dfrac{x^{k-3}}{2}\biggl(p+(p(k-3)x^2-2x-p)u_0(x)\biggr) J_{k-2}(p)\\[8pt]
&&-\;\dfrac{x^{k-4}}{2}\biggl(x+p-(x^3(k-3)-px^2(k-1)+x+p)u_0(x)\biggr)J_{k-3}(p),\\
\end{eqnarray*}
where $u_0(x)=1/{\sqrt{1-2(k-1)x^2+(k-3)^2x^4}}$.
\end{theorem}

The proof is given in Section~3. Theorem~\ref{thmain} with $k=3$ and $p=1$ shows the
generating function for the number $123$-avoiding involutions of length $n$ to be
$\frac{2x}{2x-1+\sqrt{1-4x^2}}=\sum_{n\geq0}\binom{n}{\lfloor n/2\rfloor}x^n$
(see~\cite{SS}). Also, Theorem~\ref{thmain} with $k=3$ and $p=0$ gives the number of
$123$-avoiding involutions of length $2n$ without fixed points to be
$\frac{1}{2}\binom{n}{\lfloor n/2\rfloor}$. Moreover, Theorem~\ref{thmain} with $k=4$ and
$p=1$ gives the generating function for the number $\{1234,1243\}$-avoiding involutions
of length $n$ to be $\frac{1-x}{2}+\frac{1+x}{2}\sqrt{\frac{1+2x-x^2}{1-2x-x^2}}$. In
Section~2 we present a bijection between symmetric Schr\"oder paths of length $n-1$ and
$\{1234,1243\}$-avoiding involutions of length $n$, thereby providing a combinatorial
proof of the above result with $k=4$.

\section{Symmetric Schr\"{o}der paths and $\{1234,1243\}$-avoiding involutions}

A {\it{Schr\"oder path}} of length $2n$ is a lattice path from (0,0) to $(2n,0)$
consisting of double horizontal steps $\hh=(2,0)$, up steps $\uu=(1,1)$ and down steps
$\dd=(1,-1)$ that never goes below the $x$-axis. The set of all Schr\"{o}der paths of
length $2n$ is enumerated by the $n$-th Schr\"{o}der number. Kremer~\cite[Corollary
9]{Kr} showed that $\SS_n(1243, 2143)$ is also enumerated by the $n$-th Schr\"oder number.

A Schr\"oder path of length $2n$ is called {\it{symmetric}} if it is symmetric about the
line $x=n$. Let $\paths_n$ be the collection of all such paths. In this section we give a
bijection between {\it{symmetric Schr\"{o}der paths}} of length $2n$ and the class of
involutions in $\s_{n+1}$ that avoid the patterns $\{1234,1243\}$.

We will now describe a map $\phi: \paths_n \to \I_{n+1}(1234,1243)$. 
Alongside this description will be an example of the map acting on the path 
$p=\uu\hh\uu\dd\uu\uu\dd\hh\uu\dd\dd\uu\dd\hh\dd \in \paths_9$. 
Several points may at first appear extraneous, however these will be required in the
proof of the Theorem which follows.

\centerline{\scalebox{0.6}{\includegraphics{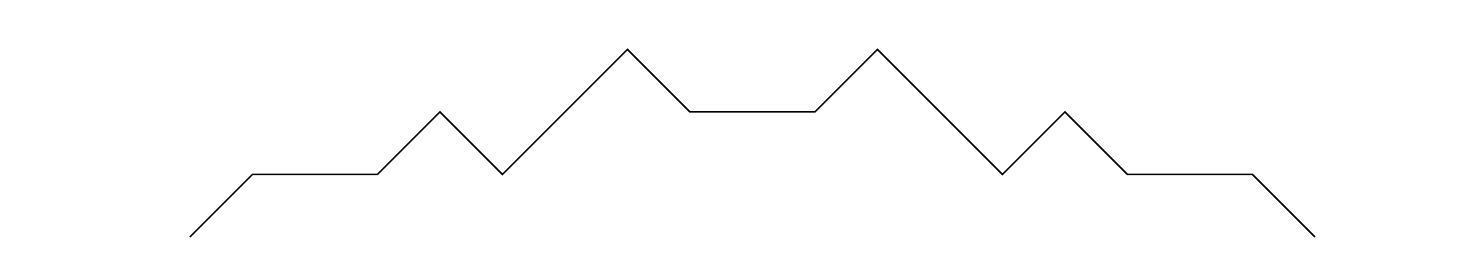}}} 

Given
$p\in\paths_n$, let $p'$ be the word of length $n+1$ obtained from $p$ by first appending
a horizontal step $\hh$ to the end, then replacing all occurrences of $\uu\dd$ in $p$ by
$\rr$, and finally deleting all remaining $\dd$'s. 
This is equivalent to projecting the steps of the path onto the diagonal, and replacing 
any $\uu$'s that are followed by a $\dd$ by $\rr$.

\centerline{\scalebox{0.6}{\includegraphics{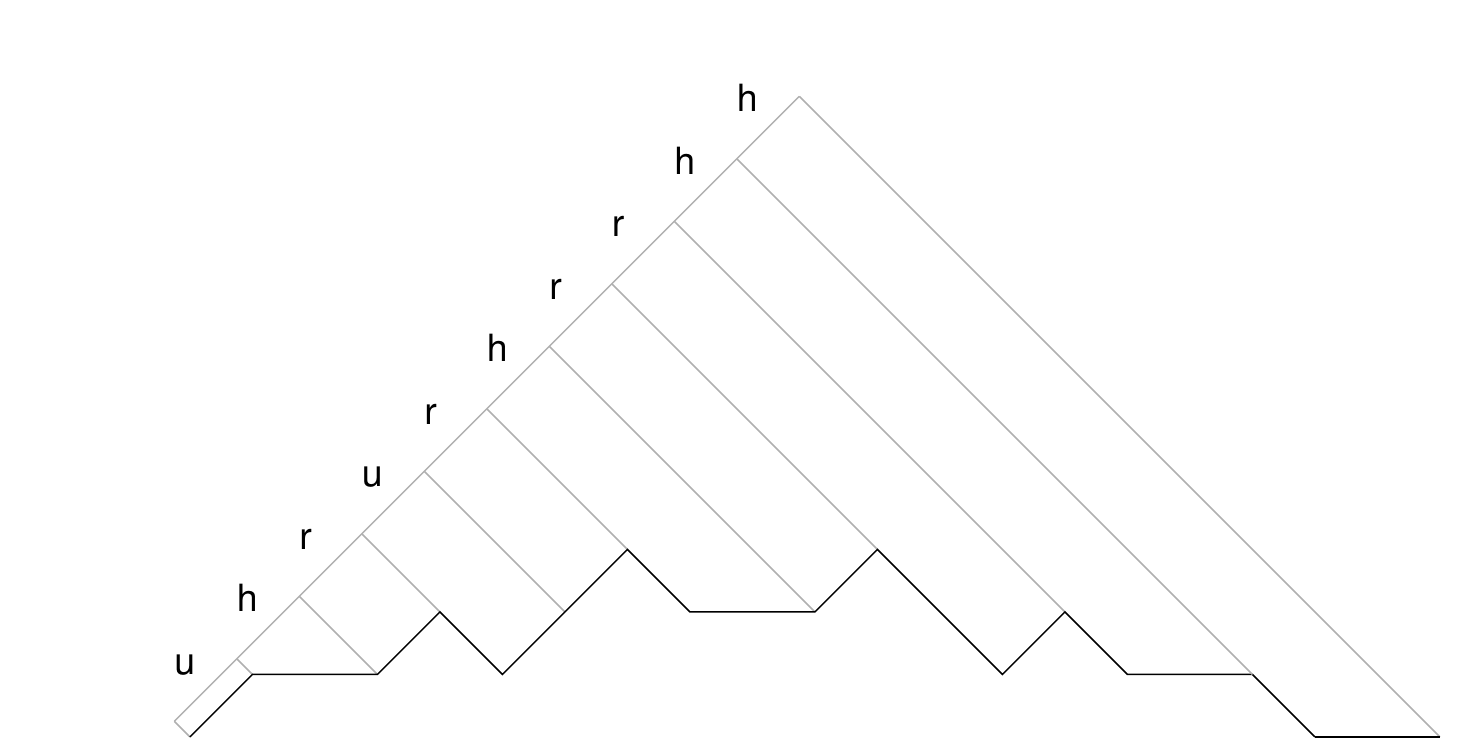}}} 

From the diagram the example path is $p'\,=\, \uu\hh\rr\uu\rr\hh\rr\rr\hh\hh$.
In general, $p'=p_1\ldots p_{n+1} \in \{\uu,\rr,\hh\}^{n+1}$. 
Now form the three sets $A_{\hh},A_{\rr},A_{\uu}$ according to
the rule $i \in A_{x}$ if $p_i = x$. 

\centerline{\scalebox{0.6}{\includegraphics{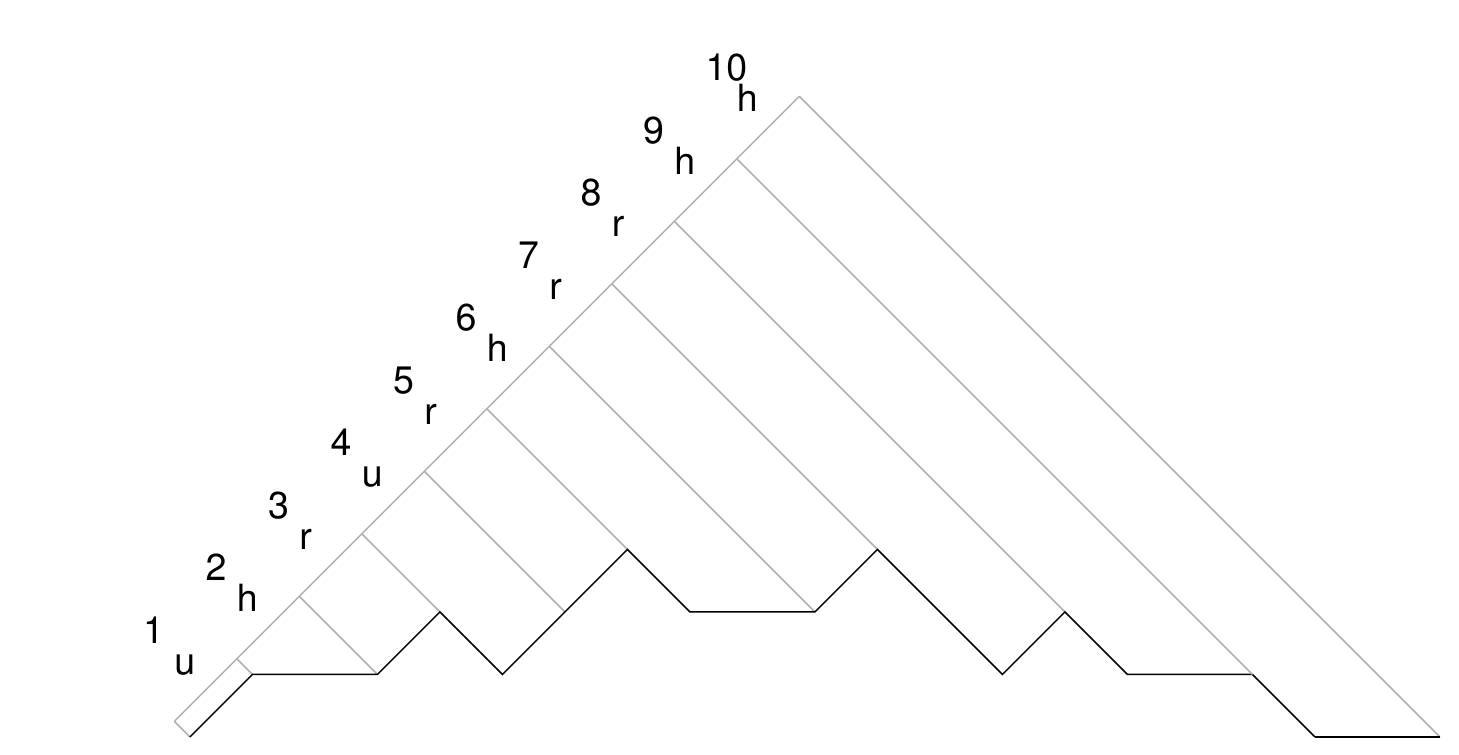}}} 

In the example, $A_{\hh}\,=\, \{2,6,9,10\}$, $A_{\rr}\,=\, \{3,5,7,8\}$ and $A_{\uu} \, = \, \{1,4\}$. 

Starting with the largest entry in $A_{\hh}$,
replace all $\hh$'s in $p'$ (from left to right) with the entries from $A_{\hh}$ so that
this sequence is decreasing. 
This is equivalent to forming a sequence of transpositions (or a fixed point), 
the first of which will have as first entry the index of the first $\hh$ in $p'$
and whose second entry is the index of the last $\hh$ in $p'$.
The second transposition (or fixed point) will have as first entry the index of the second $\hh$ in $p'$
and whose second entry is the index of the second last $\hh$ in $p'$, and so forth.

Since $A_{\hh}\,=\, \{2,6,9,10\}$ in the example, we have the transpositions $(2,10)$ and $(6,9)$.

\centerline{\scalebox{0.6}{\includegraphics{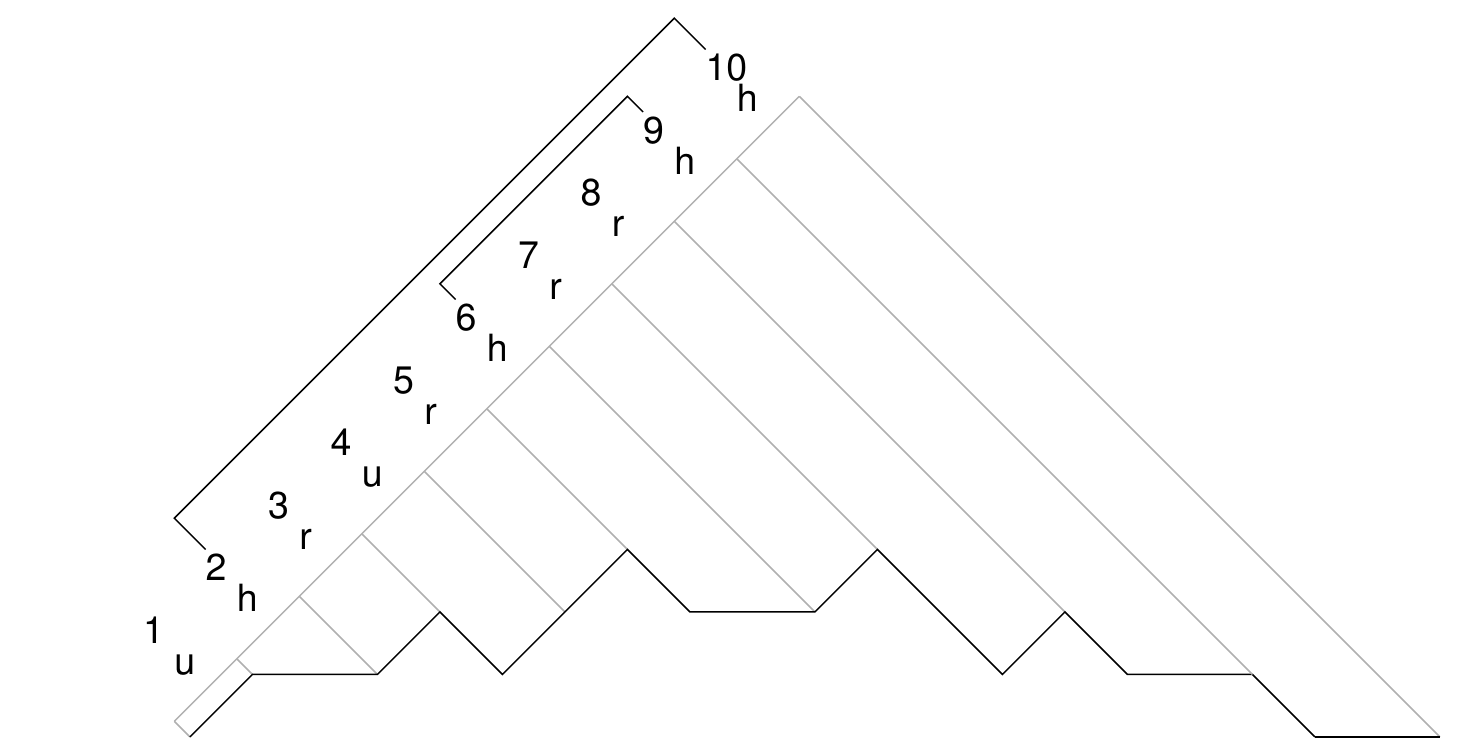}}} 

Do likewise for the sets $A_{\rr}$ and $A_{\uu}$. The
{\it{label}} of a particular $\hh$, $\rr$ or $\uu$ in $p'$ is the value that replaces it.
Call the resulting permutation $\pi = \phi(p)$, i.e. the labels of $p'$ read from left to
right. 

The label of the first $\hh$ in $p'$ above is 10 and the label of the third $\hh$ in
$p'$ above is 6.
In the example, since $A_{\hh}\,=\, \{2,6,9,10\}$, $A_{\rr}\,=\, \{3,5,7,8\}$ and $A_{\uu} \, = \, \{1,4\}$,
$\pi=\phi(p)$ is the permutation (involution) with cycles $(2,10), (6,9), (3,8), (5,7)$ and $(1,4)$.
Thus we have $\phi(p) \,=\, (4,10,8,1,7,9,5,3,6,2).$

\centerline{\scalebox{0.6}{\includegraphics{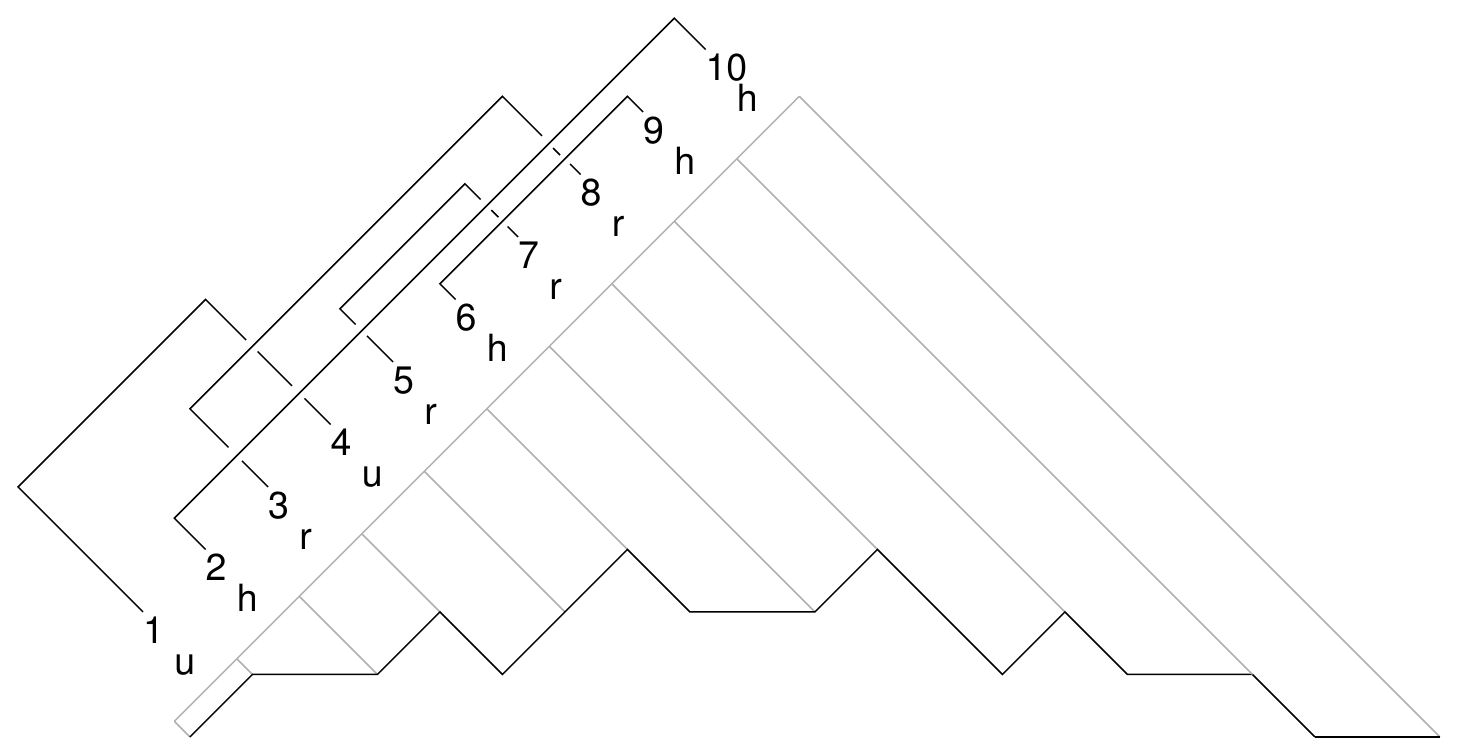}}} 

Another way to see this construction is in terms of layers of right-to-left maxima.
An element $\pi_i$ of a permutation $\pi$ is called a {\it right-to-left maximum} 
if it is greater than all elements that follow it, i.e. $\pi_i>\pi_j$ for all $j>i$.
We define successively the $r$-right-to-left maxima for a permutation $\pi\in\SS_n$. 
Let $\pi^{(1)}$ be the word consisting of all elements of $\pi$. 
For $r\ge 1$, the right-to-left maxima of $\pi^{(r)}$ are called {\it $r$-right-to-left maxima} of $\pi$.
Let $\pi^{(r+1)}$ be the subword obtained from $\pi^{(r)}$ by removing all $r$-right-to-left maxima.
For example, the permutation $\pi=  674583912\in \SS_{9}$ has the 1-right-to-left maxima 9 and 2; 
the 2-right-to-left maxima 8, 3 and 1; the 3-right-to-left maxima 7 and 5; and the 4-right-to-left maxima 6 and 4.
Note that the $r$-right-to-left maxima of $\pi$ form a decreasing subsequence for each $r$.

Let $\pi$ be the unique permutation (in fact, it will be an involution) in
$\SS_{n+1}$ with 1 right-to-left maxima $A_{\hh}$, 2 right-to-left maxima $A_{\rr}$ and 3
right-to-left maxima $A_{\uu}$. If $A_{\hh}=\{x_1,\ldots , x_{\alpha}\}$, $A_{\rr}=\{y_1,\ldots , y_{\beta}\}$
and $A_{\uu}=\{z_1,\ldots , z_{\gamma}\}$ then the cycles
of $\pi$ are 
$$(x_1\, x_{\alpha})\, (x_2 \, x_{\alpha -1})\, \cdots\, (y_1 \, y_{\beta}) \, (y_2 \, y_{\beta -1})\, \cdots \,(z_1 \, z_{\gamma}) \, (z_2\, z_{\gamma -1}) \, \cdots.$$
We point out that if $\gamma = 2m+1$ then $(x_{m+1})$ will be a fixed point. Consequently there will be at most three fixed points in the resulting involution.

The inverse map $\phi^{-1}$ is described by means of an example.
Consider $$\pi\, =\,(10,5,8,7,2,6,4,3,9,1)\in\I_{10}(1234,1243).$$ The sets of
1, 2 and 3 right-to-left maxima are $A_{\hh}=\{1,9,10\}$, $A_{\rr}=\{3,4,6,7,8\}$ and $A_{\uu}=\{2,5\}$, respectively.
This gives $p'= \hh\uu\rr\rr\uu\rr\rr\rr\hh\hh$. 
After removing the final $\hh$, we have

\centerline{\scalebox{0.6}{\includegraphics{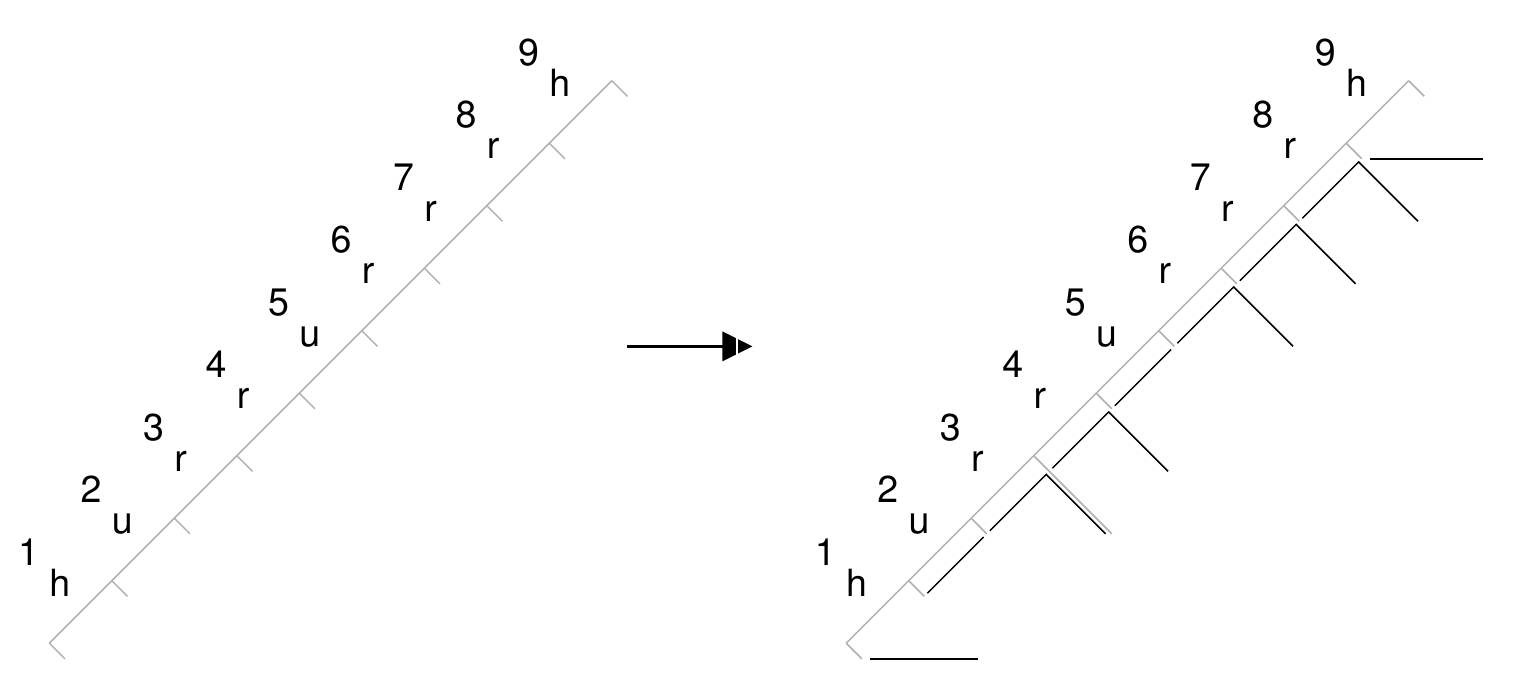}}} 

Beginning with the last step (at position 9), we push this down so that it is symmetric with the first entry.
We then move the second step of $p'$ down to meet the path. It is $\uu$ so there must be a $\dd$ inserted 
at the opposite end so the path is symmetric.

\centerline{\scalebox{0.6}{\includegraphics{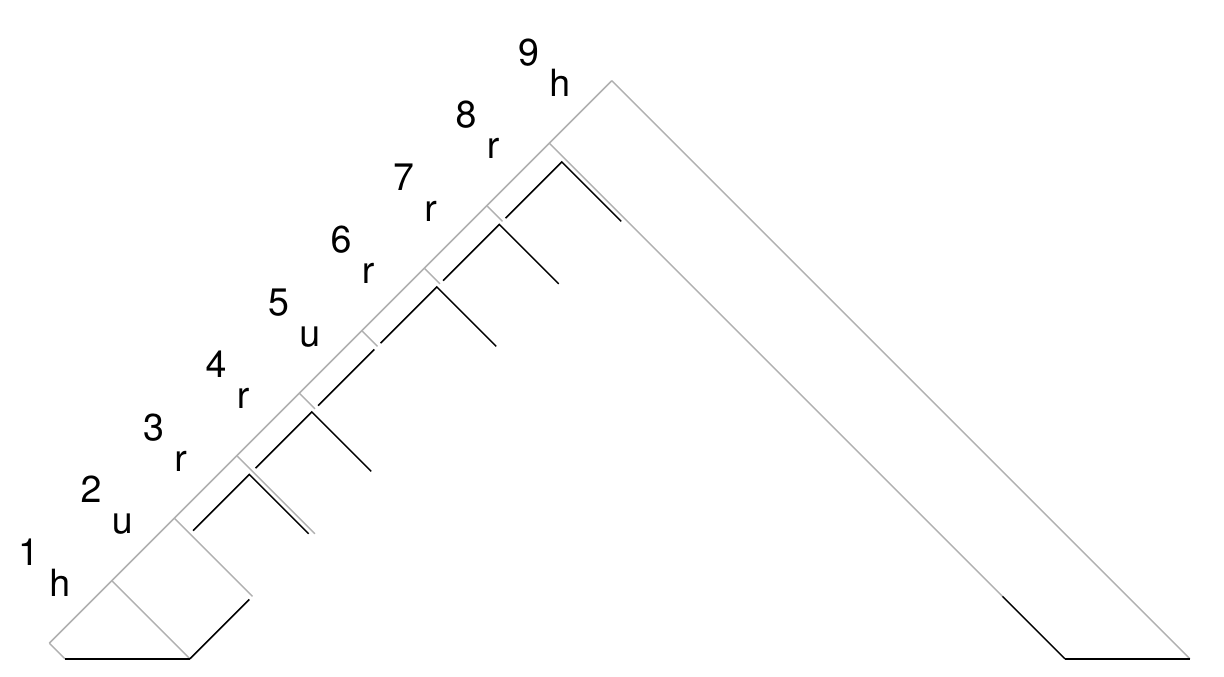}}} 

Next we move the $\rr$ at position 3 down to touch the evolving path, and move the $\rr$ at position 8 down to 
meet the path above the $\dd$ step.

\centerline{\scalebox{0.6}{\includegraphics{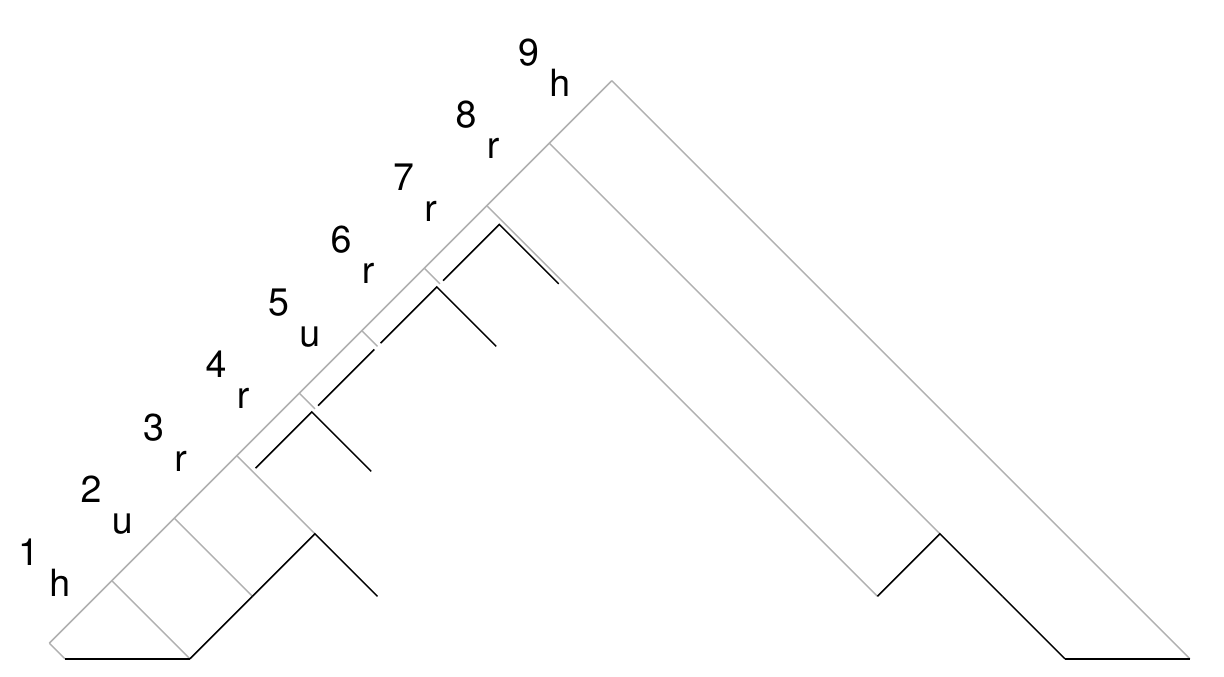}}} 

The $\uu$ at position 4 is moved next but we must insert a $\dd$ step between positions 6 and 7 to ensure the
path is symmetric.

\centerline{\scalebox{0.6}{\includegraphics{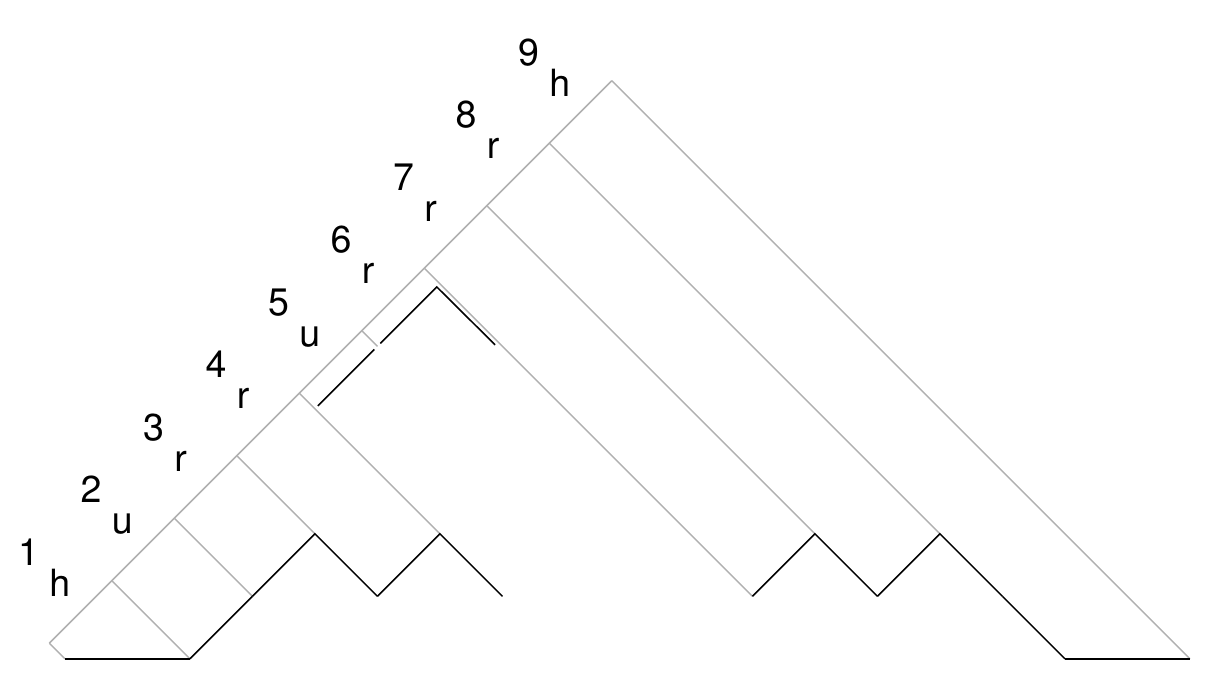}}} 

Finally, move the remaining pieces down, inserting $\dd$'s where appropriate.

\centerline{\scalebox{0.6}{\includegraphics{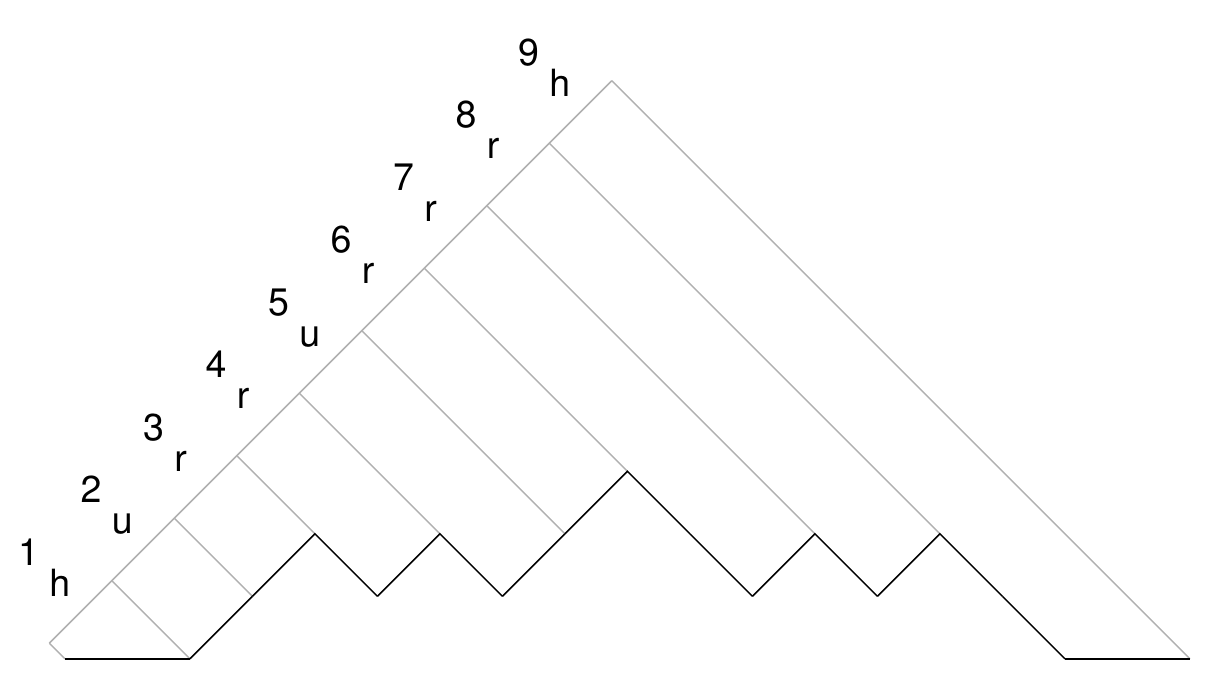}}} 

Thus we have 
$p\;=\;\phi^{-1}(\pi)\;=\;\hh\uu\uu\dd\uu\dd\uu\uu\dd\dd\uu\dd\uu\dd\dd\hh.$

\begin{theorem}
The map $\phi: \paths_n \to \I_{n+1}(1234,1243)$ is a bijection.
\end{theorem}

\begin{proof}
We first show that for any $p \in \paths_n$, the corresponding $\pi=\phi(p) \in
\I_{n+1}(1234,1243)$.

Let $p \in \paths_n$ and $p'=p_1\cdots p_k$ be the corresponding word on the alphabet
$\{\uu,\rr,\hh\}$. Suppose that $A_{\hh} = \{i_1,\ldots , i_\ell\}$. Then  it is clear
that $\pi_{i_j} = i_{\ell+1-j}$ for all $1\leq j \leq \ell$. The same is true for the
sets $A_{\rr}$ and $A_{\uu}$ so $\pi$ is an involution.

From the labelling scheme above, the resulting permutation $\pi$ has, at most, three
levels of right-to-left maxima. It is therefore 1234 avoiding. To show that $\pi$ is
1243-avoiding, suppose $\pi$ contains a 1243 pattern $\pi_i \pi_j \pi_k \pi_{\ell}$,
where $\pi_i < \pi_j < \pi_{\ell} < \pi_k$ and $i<j<k<\ell$. 
The 4 and 3 of the pattern ($\pi_k$ and $\pi_{\ell}$) must be in the 1 right-to-left maxima $A_{\hh}$.
Similarly, the 2 (resp. 1) of the pattern must be in the 2 (resp. 3) right-to-left maxima $A_{\rr}$ (resp. $A_{\uu})$.
Then $\pi_j\pi_k\pi_{\ell}$
is a 132-pattern with $\pi_j \in A_{\rr}$ and $\pi_{k},\pi_{\ell} \in A_{\hh}$. This is
not possible for the following reason: Given $\rr$ in $p'$ with label $r_1$, let
$h_1,h_2,\ldots $ be the sequence of labels of $\hh$'s to the right of $\rr$. Then
$h_1>r_1>h_2,h_3,\ldots$. This statement is easily seen by removing all $\uu$'s and the
suffix $\hh$ from $p'$ and relabelling. (This relabelling always gives a monotone
decreasing sequence.) The fact that the label of the first $\hh$ after $\rr$ is greater
than the label of the $\rr$ is due to the appended $\hh$.

We now show how to construct the unique path $p$ corresponding to $\pi \in
\I_{n+1}(1234,1243)$. For such a permutation, let $A_{\hh}$, $A_{\rr}$ and $A_{\uu}$ be
the 1, 2 and 3 right-to-left maxima of $\pi$, respectively. Insert $\hh$ at position $i$
of $p'$ if $i \in A_{\hh}$ and do likewise for the sets $A_{\rr}$ and $A_{\uu}$. Remove
the suffix $\hh$ from $p'$ (it is a suffix since $(n+1)$ is one of the 1 right-to-left
maxima). From right to left in $p'$, insert a $\dd$ where there is a corresponding $\uu$
and finish by replacing all occurrences of $\rr$ with $\uu\dd$. (As was done in the example that preceded the Theorem.)
We note that for each $p'$ there will be several Schr\"oder paths to which it may correspond, however
only one of these is symmetric.
\end{proof}

From the construction, we also have the following statistics of $\{1234, 1243\}$-avoiding
involutions:

\begin{corollary}
Let $p \in \paths_n$ with $h$ steps $\hh$, $r$ steps $\uu\dd$, and $u$ steps $\uu$ that are not directly followed by a $\dd$ step.
Let $\pi=\phi(p) \in \mathcal{I}_n(1234, 1243)$.
\begin{enumerate}
\item[(1)] The number of right-to-left maxima of $\pi$ is $h+1$.
\item[(2)] The number of 2 right-to-left maxima of $\pi$ is $r$.
\item[(2)] The number of 3 right-to-left maxima of $\pi$ is $u$.
\item[(3)] The number of fixed points of $\pi$ is $((1+h) \mbox{ mod } 2) \; +\;( r\mbox{ mod }2) \; + \; (u \mbox{ mod } 2)$.
\end{enumerate}
\end{corollary}

\begin{question}
What statistic on $\pi=\phi(p)$ corresponds to the height of the path $p$?
\end{question}

\section{Proof of Theorem~\ref{thmain}}
To present the proof of Theorem~\ref{thmain}, we must first consider the enumeration
problem for the number $\mathcal{F}_k$-avoiding involutions according to length and
number of fixed points, where $\mathcal{F}_k$ is the set of all permutations
$\sigma\in\SS_k$ with $\sigma_1=1$.

\subsection{Involutions avoiding $\mathcal{F}_k$}
In this subsection we present an explicit formula for the number of involutions that
avoid all the patterns in $\mathcal{F}_k$. To do so we require some new notation. Define
$f_k(n)$ to be the number of involutions $\pi\in\mathcal{I}_n(\mathcal{F}_k)$. Given
$t\in[n]$, we also define
 $$f_{k;m}(n;t)=\#\{\pi\in\mathcal{I}_n(\mathcal{F}_k)\mid
 \pi_1=t\mbox{ and }\pi\mbox{ contains }m\mbox{ fixed points}\}.$$
Let $f_k(n;t)=f_k(n,p;t)$ and $f_k(n)=f_k(n,p)$ be the polynomials
$\sum_{m=0}^nf_{k;m}(n;t)p^m$ and $\sum_{t=1}^nf_{k}(n;t)$, respectively. We denote by
$F_k(x,p)$ the generating function for the sequence $f_k(n,p)$, that is
$F_k(x,p)=\sum_{n\geq0} f_k(n,p)x^n$.

\begin{theorem}\label{thmm1}
We have
$$F_k(x,p)\;=\;\sum_{j=0}^{k-2}J_j(p)x^j+\frac{x^{k-1}}{1-(k-1)x^2}((k-1)J_{k-2}(p)x+J_{k-1}(p)).$$
Moreover, the number of involutions of length $k+2n$ (resp. $k+2n-1$) that avoid all the
patterns in $\mathcal{F}_k$ is given by $(k-1)^{n+1}I_{k-2}$ (resp. $(k-1)^{n}I_{k-1}$),
for all $n\geq0$.
\end{theorem}

\begin{proof}
Let $\pi\in\SS_n$ be a permutation that avoids all patterns in $\mathcal{F}_k$. We have
$\pi_1\geq n+2-k$. Thus $\pi\in\mathcal{I}_n(\mathcal{F}_k)$ with $\pi_1=t\geq n+2-k$ if
and only if $\pi_2\ldots\pi_{t-1}\pi_{t+1}\ldots\pi_n$ is an involution on the numbers
$2,\ldots,t-1,t+1,\ldots,n$ that avoids all the patterns in $\mathcal{F}_k$. Hence,
$f_k(n;j)=f_k(n-2)$ for all $j=n+2-k,n+3-k,\ldots,n$, and $f_k(n,j)=0$ for all
$j=1,2,\ldots,n+1-k$, where $n\geq k$. Thus, for $n\geq k$,
 $$f_k(n)\;=\;(k-1)f_k(n-2).$$
Using the initial conditions $f_k(j)=J_{j}(p)$, $j=1,2,\ldots,k-1$, we find that
$f_k(k+2j)=(k-1)^{j+1}J_{k-1}(p)$ and $f_k(k+2j-1)=(k-1)^jJ_{k-2}(p)$ for all $j\geq0$.
Rewriting these formulas in terms of generating functions we obtain
 $$F_k(x,p)=\sum_{j=0}^{k-2}J_j(p)x^j+\frac{x^{k-1}}{1-(k-1)x^2}((k-1)J_{k-2}(p)x+J_{k-1}(p)),$$
as claimed.
\end{proof}

\subsection{Involutions avoiding $\mathcal{A}_k$}
In this subsection we prove Theorem~\ref{thmain}. In order to do this, define $g_k(n)$ to
be the number of involutions $\pi\in\mathcal{I}_n(\mathcal{A}_k)$ and given
$t_1,t_2,\ldots,t_m\in\mathbb{N}$, we also define
 $$g_k(n;t_1,t_2,\ldots,t_m) \;=\; \#\{\pi_1\ldots\pi_n\in\mathcal{I}_n(\mathcal{A}_k)\mid \pi_1\ldots\pi_m
 =t_1\ldots t_m\}.$$
\begin{lemma}\label{lem1}
Let $k\geq3$. For all $3\leq t\leq n+1-k$,
 $$g_k(n;t)\;=\;(k-2)g_k(n-2;t-1)+\sum_{j=1}^{t-2}g_k(n-2;j),$$
with the initial conditions $g_k(n;1)=f_{k-1}(n-1)$, $g_k(n;2)=f_{k-1}(n-2)$, and
$g_k(n;t)=g_k(n-2)$ for all $t=n+2-k,n+3-k,\ldots,n$.
\end{lemma}

\begin{proof}
Let $\pi$ be any involution of length $n$ that avoids all patterns in $\mathcal{A}_k$
with $\pi_1=t$. Now let us consider all possible values of $t$. If $t=1$ then
$\pi\in\mathcal{I}_n(\mathcal{A}_k)$ if and only if
$(\pi_2-1)(\pi_3-1)\ldots(\pi_n-1)\in\mathcal{I}_{n-1}(\mathcal{F}_{k-1})$. If $t=2$ then
$\pi\in\mathcal{I}_n(\mathcal{A}_k)$ if and only if
$(\pi_3-2)(\pi_4-2)\ldots(\pi_n-2)\in\mathcal{I}_{n-2}(\mathcal{F}_{k-1})$. Now assume
that $3\leq t\leq n+1-k$, then from the above definitions
$$\begin{array}{l}
g_k(n;t) \;=\; g_k(n;t,1)+\ldots+g_k(n;t,t-1)+g_k(n;t,t+1)+\cdots+g_k(n;t,n).
\end{array}$$
But any involution $\pi$ satisfying $\pi_1<\pi_2\leq n+2-k$ contains a pattern from the
set $\mathcal{A}_k$ (see the subsequence of the letters $\pi_1,\pi_2,
n+3-k,n+4-k,\ldots,n$ in $\pi$). Thus $g_k(n;t,r)=0$ for all $t<r\leq n+2-k$ and so
$$\begin{array}{l}
g_k(n;t) \;=\; g_k(n;t,1)+\ldots+g_k(n;t,t-1)+g_k(n;t,n+3-k)+\cdots+g_k(n;t,n).
\end{array}$$
Also, if $\pi$ is an involution in $\mathcal{I}_n$ with $\pi_1=t$ and $\pi_2=r\geq
n+3-k$, then the entry $r$ does not appear in any occurrence of $\tau\in\mathcal{A}_k$ in
$\pi$. Thus, there exists a bijection between the set of involutions
$\pi\in\mathcal{I}_n(\mathcal{A}_k)$ with $\pi_1=t$ and $\pi_2=r\geq n+3-k$ and the set
of involutions $\pi'\in\mathcal{I}_{n-2}(\mathcal{A}_k)$ with $\pi'=t-1$. Therefore
$g_k(n;t,r)=g_k(n-2;t-1)$ which gives
$$\begin{array}{l}
g_k(n;t) \;=\; g_k(n;t,1)+\ldots+g_k(n;t,t-1)+(k-2)g_k(n-2;t-1).
\end{array}$$
Also, if $\pi$ is an involution in $\mathcal{I}_n$ with $\pi_1=t$, $\pi_2=r<t$ and if
$ta_2\ldots a_k$ is an occurrence of a pattern from the set $\mathcal{A}_k$ in $\pi$,
then $ra_2\ldots a_k$ is an occurrence of a pattern from the set $\mathcal{A}_k$ in
$\pi$. Thus, there exists a bijection between the set of involutions
$\pi\in\mathcal{I}_n(\mathcal{A}_k)$ with $\pi_1=t$ and $\pi_2=r<t$ and the set of
involutions $\pi'\in\mathcal{I}_{n-2}(\mathcal{A}_k)$ with $\pi'_1=r-1$. Therefore
$g_k(n;t,r)=g_k(n-2;r-1)$ which gives
$$\begin{array}{l}
g_k(n;t)\;=\;(k-2)g_k(n-2;t-1)+\displaystyle\sum_{j=1}^{t-2} g_k(n-2;j),
\end{array}$$
as required. Finally, if $\pi$ is an involution in $\mathcal{I}_n$ with $\pi_1=t\geq
n+2-k$, then the entry $t$ does not appear in any occurrence of $\tau\in\mathcal{A}_k$ in
$\pi$. Thus, there exists a bijection between the set of involutions
$\pi\in\mathcal{I}_n(\mathcal{A}_k)$ with $\pi_1=t\geq n+2-k$ and the set of involutions
$\pi'\in\mathcal{I}_{n-2}(\mathcal{A}_k)$. Therefore $g_k(n;t)=g_k(n-2)$, as claimed.
\end{proof}

Let $G_k(n;v)$ be the polynomial $\sum_{t=1}^ng_k(n;t)v^{t-1}$. Rewriting the above lemma
in terms of the polynomials $G_k(n;v)$ we have the following recurrence relation.
\begin{lemma}\label{lem2}
Let $k\geq3$. For all $n\geq k$,
$$\begin{array}{lcl}
\lefteqn{G_k(n;v)}\\
&=&f_{k-1}(n-1)+vf_{k-1}(n-2)-v(k-2)f_{k-1}(n-3)+\left(\frac{v^2}{1-v}+(k-2)v\right)G_k(n-2;v)\\
&&-\frac{v^n}{1-v}G_k(n-2;1)+\frac{v^{n-1}}{1-v}\left(k-2+\frac{v-v^{3-k}}{1-v}\right)G_k(n-4;1),
\end{array}$$
where $G_k(n;v)=I_{n-1}+\frac{v-v^n}{1-v}I_{n-2}$ for all $n=0,1,\ldots,k-1$.
\end{lemma}
\begin{proof}

Lemma~\ref{lem1} gives
$$\begin{array}{cl}
\lefteqn{G_k(n;v)}\\
=&f_{k-1}(n-1)+vf_{k-1}(n-2)+\displaystyle\sum_{t=2}^{n-k}v^t\left((k-2)G_k(n-2;t)+\displaystyle\sum_{j=1}^{t-1}G_k(n-2;j)\right)\\
&+ G_k(n-2;1) \displaystyle\sum_{t=n+1-k}^{n-1}v^t\\
=& f_{k-1}(n-1)+vf_{k-1}(n-2)+\frac{v^2}{1-v}\left(G_k(n-2;v)-G_k(n-4;1)\displaystyle\sum_{j=n-1-k}^{n-3}v^j\right)\\
&-\frac{v^{n+1-k}}{1-v}(G_k(n-2;1)-(k-1)G_k(n-4;1))+\frac{v^{n+1-k}-v^n}{1-v}G_k(n-2;1)\\
&+(k-2)v\left(G_k(n-2;v)-f_{k-1}(n-3)-G_k(n-4;1)\displaystyle\sum_{j=n-k}^{n-3}v^j\right),
\end{array}$$
which is equivalent to
$$\begin{array}{ll}
\lefteqn{G_k(n;v)}\\
=&f_{k-1}(n-1)+vf_{k-1}(n-2)-v(k-2)f_{k-1}(n-3)+\left(\frac{v^2}{1-v}+(k-2)v\right)G_k(n-2;v)\\
&-\frac{v^n}{1-v}G_k(n-2;1)+\frac{v^{n-1}}{1-v}\left(k-2+\frac{v-v^{3-k}}{1-v}\right)G_k(n-4;1).
\end{array}$$
To find the value of $G_k(n;v)$ for $n\leq k-1$, let $\pi$ be any involution with
$\pi_1=t$. If $t=1$ then there are $I_{n-1}$ involutions, whereas if $t>1$ there are
$I_{n-2}$ involutions, hence
$G_k(n;v)=v^0I_{n-1}+\sum_{t=2}^nv^{t-1}I_{n-2}=I_{n-1}+\frac{v-v^n}{1-v}I_{n-2}$, as
required.
\end{proof}

Lemma~\ref{lem2} can be generalised as follows; let $g_{k;m}(n;t)$ be the number of
involutions $\pi\in\mathcal{I}_n(\mathcal{A}_k)$ such that $\pi_1=t$ and $\pi$ contains
exactly $m$ fixed points. Define $G_k(n;t;p)=\sum_{m=0}^ng_{k;m}(n;t)p^m$ and
$G_k(n;v,p)=\sum_{t=1}^nG_k(n;t;p)v^{t-1}$. Using the same arguments as those in the
proofs of Lemma~\ref{lem1} and Lemma~\ref{lem2}, while carefully considering the number
of fixed points, we have the following result.

\begin{lemma}\label{lem3}
Let $k\geq3$. For all $n\geq k$,
$$\begin{array}{ll}
\lefteqn{G_k(n;v,p)}\\
=&pf_{k-1}(n-1)+vf_{k-1}(n-2)-pv(k-2)f_{k-1}(n-3)+\left(\frac{v^2}{1-v}+(k-2)v\right)G_k(n-2;v,p)\\
&-\frac{v^n}{1-v}G_k(n-2;1,p)+\frac{v^{n-1}}{1-v}\left(k-2+\frac{v-v^{3-k}}{1-v}\right)G_k(n-4;1,p),
\end{array}$$
where $G_k(n;v,p)=pJ_{n-1}(p)+\frac{v-v^n}{1-v}J_{n-2}(p)$ for all $n=0,1,\ldots,k-1$.
\end{lemma}

Let $G_k(x,v,p)=\sum_{n\geq0}G_k(n;v,p)x^n$ be the generating function for the sequence
$G_k(n;v,p)$. Define $J_i(v,p)$ to be the polynomial $\sum d_{tr}v^tp^r$ where $d_{tr}$
is the number of involutions $\pi\in\mathcal{I}_i$ such that $\pi_1=t+1$ and $\pi$
contains exactly $r$ fixed points. Rewriting the recurrence relation in the statement of
Lemma~\ref{lem3} in terms of generating functions we obtain
$$\begin{array}{l}
\lefteqn{G_k(x,v,p)\;=\;}\\
\sum\limits_{j=0}^{k-1}J_j(v,p)x^j+px\left(F_{k-1}(x,p)-\sum\limits_{j=0}^{k-2}J_j(p)x^j\right)+vx^2\left(F_{k-1}(x,p)-\sum\limits_{j=0}^{k-3}J_j(p)x^j\right)\\
-(k-2)pvx^3\left(F_{k-1}(x,p)-\sum\limits_{j=0}^{k-4}J_j(p)x^j\right)-\frac{v^2x^2}{1-v}\left(G_k(xv,1,p)-\sum\limits_{j=0}^{k-3}J_j(p)(xv)^j\right)\\
+vx^2\left(\frac{v}{1-v}+k-2\right)\left(G_k(x,v,p)-\sum\limits_{j=0}^{k-3}J_j(v,p)x^j\right)\\
+\frac{(k-2)v^3x^4}{1-v}\left(G_k(xv,1,p)-\sum\limits_{j=0}^{k-5}J_j(p)(xv)^j\right)
-\frac{x^4(1-v^{k-2})}{v^{k-6}(1-v)^2}\left(G_k(xv,1,p)-\sum\limits_{j=0}^{k-5}J_j(p)(xv)^j\right),
\end{array}$$
which is equivalent to
$$\begin{array}{l}
\left(1-\frac{x^2}{1-v}-(k-2)\frac{x^2}{v}\right)G_k(x/v,v,p) \; = \; \\
-\frac{x^2}{1-v}\left(1-(k-2)\frac{x^2}{v}+\frac{x^2(1-v^{k-2})}{v^{k-2}(1-v)}\right)G_k(x,1,p)\\
+\sum\limits_{j=0}^{k-1}J_j(v,p)\frac{x^j}{v^j}+\frac{px}{v}\left(F_{k-1}(x/v,p)-\sum\limits_{j=0}^{k-2}J_j(p)\frac{x^j}{v^j}\right)+\frac{x^2}{v}\left(F_{k-1}(x/v,p)-\sum\limits_{j=0}^{k-3}J_j(p)\frac{x^j}{v^j}\right)\\
-(k-2)p\frac{x^3}{v^2}\left(F_{k-1}(x/v,p)-\sum\limits_{j=0}^{k-4}J_j(p)\frac{x^j}{v^j}\right)+\frac{x^2}{1-v}\sum\limits_{j=0}^{k-3}J_j(p)x^j\\
-\frac{x^2}{v}\left(\frac{v}{1-v}+k-2\right)\sum\limits_{j=0}^{k-3}J_j(v,p)\frac{x^j}{v^j}
-\frac{(k-2)x^4}{v(1-v)}\sum\limits_{j=0}^{k-5}J_j(p)x^j+\frac{x^4(1-v^{k-2})}{v^{k-2}(1-v)^2}\sum\limits_{j=0}^{k-5}J_j(p)x^j.
\end{array}$$
To solve this functional equation, we substitute
$$v:=v_0=\frac{1}{2}\left(1+(k-3)x^2+\sqrt{1-2(k-1)x^2+(k-3)^2x^4}\right),$$
where $v_0$ is the root of the coefficient of $G_k(x/v,v,p)$ above, into the above
functional equation, that is, $1-\frac{x^2}{1-v_0}-(k-2)\frac{x^2}{v_0}=0$. Since
$J_j(v,p)=pJ_{j-1}(p)+\frac{v-v^j}{1-v}J_{j-2}(p)$ for all $j=1,2,\ldots,k-1$ and
$J_0(v,p)=1$, it is routine to show (via some rather tedious algebraic manipulation) that
we obtain Theorem~\ref{thmain}.


\end{document}